\documentclass[12pt]{paper}

\usepackage{amssymb,amsfonts,amsthm,amsmath}

\usepackage{latexsym,multicol,graphicx}

\usepackage{graphicx, verbatim}
\usepackage{xcolor}
\usepackage{verbatim}

\usepackage[a4paper, hmargin=3cm, vmargin={4cm, 4cm}]{geometry}

\usepackage{fancyhdr}
\newtheorem{theorem}{Theorem}
\newtheorem*{theorem*}{Theorem}
\newtheorem{lemma}[theorem]{Lemma}
\newtheorem{question}[theorem]{Question}

\newtheorem{corollary}[theorem]{Corollary}

\renewcommand{\phi}{\varphi}

\newcommand{\eps}{\varepsilon}

\newcommand{\la}{\lambda}

\newcommand{\cP}{\mathcal P}

\newcommand{\bE}{\mathbb E}

\renewcommand{\le}{\leqslant}
\renewcommand{\ge}{\geqslant}

\newcommand{\bR}{\mathbb R}
\newcommand{\bC}{\mathbb C}
\newcommand{\bZ}{\mathbb Z}
\newcommand{\bD}{\mathbb D}
\newcommand{\bT}{\mathbb T}
\newcommand{\bN}{\mathbb N}

\newcommand{\bP}{\mathbb P}

\begin{document}

\title{Notes on the Szeg\H{o} minimum problem. \\
I. Measures with deep zeroes}

\author{Alexander Borichev
\thanks{Supported by a joint grant of Russian Foundation for Basic Research
and CNRS (projects 17-51-150005-NCNI-a and PRC CNRS/RFBR 2017-2019) and by the project ANR-18-CE40-0035.} \\
\and
Anna Kononova
\thanks{Supported by a joint grant of Russian Foundation for Basic Research
and CNRS (projects 17-51-150005-NCNI-a and PRC CNRS/RFBR 2017-2019).} \\
\and
Mikhail Sodin
\thanks{Supported by ERC
Advanced Grant 692616 and ISF Grant 382/15.} }

\maketitle


\begin{abstract}
The classical Szeg\H{o} polynomial approximation theorem states that
the polynomials are dense in the space $L^2(\rho)$, where $\rho$ is a
measure on the unit circle, if and only if the logarithmic integral of
the measure $\rho$ diverges. In this note we give a quantitative version
of Szeg\H{o}'s theorem in the special case when the divergence of the
logarithmic integral is caused by deep zeroes of the measure $\rho$ on a
sufficiently rare subset of the circle.
\end{abstract}

\section{Introduction}

Denote by $\cP$ the linear space of algebraic polynomials, and by $\cP_n$ its subspace
of polynomials of degree $n$.
Given a finite positive measure $\rho$ on the unit circle $\bT$, put
\[
e_n(\rho) = \min_{q_0, \ldots , q_{n-1}}\, \sqrt{\int_{\bT}
\bigl| q_0+q_1t + \ldots +q_{n-1}t^{n-1} +t^n \bigr|^2\, {\rm d}\rho(t)}
= \text{dist}_{L^2(\rho)}(t^n, \cP_{n-1}).
\]
Then
\[
\lim_{n\to\infty} e_n(\rho) =
\exp\Bigl( \frac12\int_\bT \log\rho'\, {\rm d}m \Bigr),
\]
where $m$ is the Lebesgue measure on $\bT$ normalized by condition $m(\bT)=1$,
and $\rho' = {\rm d}\rho/{\rm d}m$ is the Radon--Nikodym derivative. This is a classical
result, first, proven by Szeg\H{o} for absolutely continuous measures $\rho$,
and then, independently, by Verblunsky and Kolmogorov in the general case~\cite[Section~3.1]{GS}
and~\cite[Chapters~1 and 2]{Simon}. Noting that for $j\ge 0$,
$e_{n+j}(\rho)$ coincides with the distance in $L^2(\rho)$ from
$ t^{-j}$ to the linear span of $\{ t^k\colon -j+1 \le k \le n \}$,
and recalling that
the trigonometric polynomials are dense in $L^2(\rho)$, one sees
that the density of algebraic polynomials $\cP$ in $L^2(\rho)$ is equivalent to the condition
$ \displaystyle\lim_{n\to\infty} e_n(\rho) = 0$, and therefore,
to the divergence of the logarithmic integral
\[
\int_\bT \log\rho'\, {\rm d}m = -\infty\,.
\]

In these notes we will be occupied by the following question:

\begin{question}\label{quest} Suppose $\rho$ is a measure on $\bT$ with divergent logarithmic integral. Estimate the rate of decay of the sequence $e_n(\rho)$.
\end{question}

\smallskip\noindent Our interest to this question came from the linear prediction
for stationary processes.
If $\xi\colon \bZ\to\bC$ is a stationary random sequence with spectral measure $\rho$, then, according to Kolmogorov and Wiener, $e_n(\rho)$
is the error of the best mean-quadratic linear prediction of
$\xi(n)$ by $\xi(0), \ldots , \xi(n-1)$; i.e.,
\[
e_n^2(\rho) = \min_{q_0, \ldots , q_{n-1}} \bE \Bigl[ \Bigl|  \xi(n) - \sum_{0\le j \le n-1} q_j \xi(j) \Bigr|^2  \Bigr].
\]

In the case when the logarithmic integral converges, $e_n(\rho)$ has a positive limit $e_\infty(\rho)$, and dependence of
the rate of convergence on the smoothness of the density of $\rho$ is
well-understood~\cite{Gol, Ibr}. In the case of divergent logarithmic
integral the situation is quite different and not much is known. If the closed support of $\rho$ is not the whole circle, then
it is not difficult to show that $e_n(\rho)$ tends to zero at least exponentially. In the other direction, a version of the classical result
of Erd\H{o}s and Tur\'an says if $\rho' >0$ $m$-a.e. on $\bT$, then the measure
$\rho$ is regular, i.e., $e_n(\rho)^{1/n}\to 1$. Later,
stronger criteria for regularity of $\rho$ were found by Widom, Ullman, and
Stahl and Totik, see~\cite[Chapter~4]{ST}.

In these notes we show that in several special but interesting situations
it is not difficult to estimate decay of the sequence $e_n(\rho)$ using only simple
classical tools. Here, we consider the case when the divergence of the logarithmic integral is caused by deep zeroes of the measure $\rho$ on a sufficiently
rare subset of $\bT$. The results presented in this note  extend Theorems~8 and~9 from~\cite{BSW}.

Our main idea is that in the case when the measure ${\rm d}\rho = \Phi\, {\rm d}m$
has divergent logarithmic integral (i.e.,
$\displaystyle \int_\bT \log\Phi\, {\rm d}m = -\infty $),
the value $|\log e_n(\rho)|$ can be controlled by the integral
\[
\int_\bT \min\Bigl\{ \log\Bigl(\frac1{\Phi} \Bigr), A \Bigr\}\, {\rm d}m
\]
of the cut-off of $\log\Phi^{-1}$ on an appropriate large level $A$
depending on $n$. This can be viewed as a quantitative version of the regularization
of the weight $\Phi$ by $\Phi_\eps = \Phi+\eps$ with
$0<\eps\ll 1$ used by Szeg\H{o} in the proof of his theorem.
We succeeded to make this work only under additional regularity
assumptions on $\Phi$.

The toy example is the absolutely continuous measure
${\rm d}\rho=e^{-H}\, {\rm d}m$, where $H(e^{2\pi{\rm i}\theta})=h(\theta)$,
$h\colon \bR\to [0, +\infty]$ is a $1$-periodic even function, continuous and decreasing on $(0, \frac12]$, and such that
$\displaystyle \int_0 h(\theta)\,{\rm d}\theta = +\infty$. Then, under mild
assumptions on $h$, we obtain
\[
| \log e_n(\rho) | \simeq \int_\bT h_A\, {\rm d}m\,,
\]
where $h_A=\max(h, A)$, and $A=A(n)$ is a solution to the equation $n h^{-1}(a)=a$,
$h^{-1}$ is the inverse to the restriction of $h$ on $(0, \frac12]$.  
Throughout the paper we use the following notation: for positive $A$ and $B$, 
$A \lesssim B $ means that there is a positive numerical constant $C$ such that
$A \le CB $, while $A \gtrsim B $ means that $B \lesssim A $, and
$A \simeq B $ means that both $A \lesssim B $ and
$B \lesssim A $.

In the forthcoming second note, we will consider the opposite case when the bulk
of the measure $\rho$ is concentrated on a rare subset of $\bT$.

\subsubsection*{Acknowledgements}
We thank Sergei  Denisov, Fedor Nazarov, and Eero Saksman for several enlightening
discussions.

\section{Preliminaries}
Here and elsewhere, $H\colon \bT\to [0, +\infty]$ is a measurable function with
\[
\int_\bT H\, {\rm d}m = +\infty.
\]
By $\la_H (a) = m\{H>a\}$ we denote the distribution function of $H$.
For $A\ge 1$, we put $H_A(t)=\min(H(t), A)$.
To estimate from below and above $\log e_n$, we will use the integrals
\[
\int_\bT H_A\, {\rm d}m = \int_0^A \la_H (a)\, {\rm d}a
= A\la_H(A) + \int_{\{H\le A\}} H\, {\rm d}m
\]
with some $A=A(n)$.

We record several simple observations, which we frequently use throughout the paper.

\subsection{}
First, we note that under mild regularity assumptions one of the two terms
on the RHS can be discarded.  If $\la_H (a)$ satisfies
\[
\limsup_{a\to\infty} \frac{\la_H(a)}{\la_H(2a)} < 2
\]
(i.e., decays not faster than $a^{-p}$ with some $p<1$), then
\[
\int_\bT H_A\, {\rm d}m \simeq A\la_H(A).
\]

\subsection{}
On the other hand, if the function $a\mapsto a\la_H^2(a)$ does not increase
(i.e., $\la_H(a)$ decays as $1/\sqrt a$, or faster), then
\[
\int_\bT H_A\, {\rm d}m \simeq \int_{\{H\le A\}} H\, {\rm d}m,
\]
provided that $\la_H(A)$ is separated from $1$ (i.e.,
$A$ is sufficiently large).
To see this, denote by $H^*\colon [0, 1]\to [0, +\infty]$ the decreasing rearrangement
of $H$, that is, the function inverse to $\la_H$. Then, the function $s\mapsto s^2H^*(s)$
does not decrease. Letting $\alpha=\la_H(A)$ (i.e., $A=H^*(\alpha)$), we obtain
\begin{multline*}
A\la_H(A) = \alpha H^*(\alpha) \lesssim
\alpha^2 H^*(\alpha)\, \int_\alpha^1 \frac{{\rm d}s}{s^2} \\
\le \int_\alpha^1 \frac{s^2 H^*(s)}{s^2}\, {\rm d}s = \int_\alpha^1 H^*(s)\, {\rm d}s
= \int_{\{H\le A\}} H\, {\rm d}m.
\end{multline*}

\subsection{}
Furthermore, if $A/\la_H(A) \simeq B/\la_H(B)$, then
\[
\int_\bT H_A\, {\rm d}m \simeq \int_\bT H_B\, {\rm d}m.
\]
Assume, for instance, that $A\le B$ and that $B/\la_H(B) \le C\cdot A/\la_H(A)$.
Since $\la_H$ does not increase, $B\le C\cdot A$. Then,
\begin{multline*}
\int_{\{H\le B\}} H\, {\rm d}m \le \int_{\{H\le A\}} H\, {\rm d}m
+ \int_{\{A < H\le C\cdot A\}} H\, {\rm d}m \\
\le \int_{\{H\le A\}} H\, {\rm d}m + C\cdot A\la_H(A) \le C\, \int_\bT H_A\, {\rm d}m\,,
\end{multline*}
and similarly,
\[
B\la_H(B) = \frac{B}{\la_H(B)}\cdot \la_H^2(B)
\le C\cdot \frac{A}{\la_H(A)}\cdot \la_H^2(A) = C\cdot A\la_H(A).
\]
Thus,
\[
\int_\bT H_B\, {\rm d}m \le C\, \int_\bT H_A\, {\rm d}m.
\]
Since $B\ge A$, the opposite estimate is obvious.

\subsection{}
Our last remark concerns regularity of $H$. Since we will be interested only
in rather crude lower and upper bounds for $e_n(\rho)$ under conditions
${\rm d}\rho \ge e^{-H}\, {\rm d}m$ (for lower bounds) or
$\displaystyle \int_\bT e^H\, {\rm d}\rho <\infty$ (for upper bounds),
our estimates will not distinguish between the sequences $e_n(\rho)$ and
$C e_n(\rho)$, and
we can always replace the function $H$ by any function $\widetilde{H}$ with
$|\widetilde{H}-H|\le 1$ without affecting our estimates.
Keeping this in mind, we always assume that, for any positive $C$,
the equation $C\la_H(A)=A$ has a unique solution.

In Section~\ref{sec5} (Theorem~\ref{the9}) we will be using the same tacit assumption
for the equation $C\la^*_H(A)=A$, where $\la_H^*(a)$ is the length of the
longest open interval within the set $\{H>a\}$.

In the same way, we can always assume that $e_n (\rho)<1/2$.

\section{The lower bound for $e_n$ via the Remez-type inequality}

\begin{theorem}\label{thm:LB}
Suppose ${\rm d}\rho \ge e^{-H}{\rm d}m$. Then
\[
| \log e_n (\rho)| \lesssim \int_0^A \la_H(a)\,{\rm d}a,
\]
where $A=A(n)$ solves the equation $n \la_H(A)=A$.
\end{theorem}
\begin{corollary}\label{cor:LB} Suppose $H$ belongs to the weak $L^1(m)$-space, i.e.,
$\la_H(a)\lesssim 1/a$ for $a\ge 1$. Then $e_n(\rho)$ does not decay to zero faster than
a negative power of $n$.
\end{corollary}
Similarly, if $H$ belongs to the weak $L^p(m)$-space with $0<p<1$,
that is $\la_H(a) \lesssim a^{-p}$,
then
\[ | \log e_n(\rho) | \lesssim n^{\frac{1-p}{1+p}}\,.\]
If $\la_H(a) \lesssim (\log a)^{-1}$ for $a\ge e$ (in particular, if $\log H$ is
integrable), then \[ | \log e_n(\rho) | \lesssim \frac{n}{\log^2 n}\,,\]
and so on, until we arrive at the classical Erd\H{o}s-Tur\'an theorem, which states that
\[
|\log e_n(\rho)|=o(n), \quad n\to\infty,
\]
provided that $H<+\infty$ a.e.
on $\bT$ (that is, $\la_H(a)\to 0$ as $a\to\infty$).

\subsection{Proof of Theorem~\ref{thm:LB}}

Let $P$ be an extremal algebraic polynomial of degree $n$ such that $P(0)=1$ and
\[
\int_\bT |P|^2\, {\rm d}\rho = e_n(\rho)^2.
\]
Then
\begin{align*}
0 &\le \int_\bT \log|P|^2\, {\rm d}m \\
&=\int_{\{H>A\}} \log|P|^2\, {\rm d}m
+ \int_{\{H\le A\}} \log\bigl( |P|^2 e^{-H} \bigr)\, {\rm d}m
+ \int_{\{H\le A\}} H\, {\rm d}m.
\end{align*}
Estimating the first and the second integrals on the RHS we let $n$ be so large that
$ \int_{\{H\le A\}} H\, {\rm d}m \ge 1$ and
\begin{equation}
e^A > \rho(\bT).
\label{ko1}
\end{equation}
By Jensen's inequality,
\begin{align*}
\int_{\{H\le A\}} \log\bigl( |P|^2 e^{-H}\bigr)\, {\rm d}m
&\le \frac1e + \log\Bigl( \int_{\{H\le A\}} |P|^2 e^{-H}\, {\rm d}m \Bigr) \\
&\le \frac1e + \log\Bigl( \int_\bT |P|^2 \, {\rm d}\rho \Bigr) = \frac1e + 2 \log e_n(\rho),
\end{align*}
and similarly,
\begin{align*}
\int_{\{H>A\}} \log|P|^2\, {\rm d}m &=
\la_H(A) \Bigl( \frac1{\la_H(A)}\,\int_{\{H>A\}} \log |P|^2\, {\rm d}m \Bigr) \\
&\le \la_H(A) \log\Bigl( \frac1{\la_H(A)}\,\int_{\{H>A\}} |P|^2\, {\rm d}m \Bigr) \\
&\le \frac1e + \la_H(A)\log\Bigl( \int_\bT |P|^2\, {\rm d}m \Bigr).
\end{align*}
Next, applying the $L^2$-version of the classical Remez inequality
(which follows, for instance, from a more general Nazarov's result~\cite{Nazarov}),
we obtain
\begin{align*}
\int_\bT |P|^2\, {\rm d}m &\le e^{Cn\la_H(A)}\int_{\{H\le A\}} |P|^2\, {\rm d}m \\
&\le e^{Cn\la_H(A)+A}\int_{\{H\le A\}} |P|^2e^{-H}\, {\rm d}m \\
&\le e^{Cn\la_H(A)+A}\int_{\bT} |P|^2\, {\rm d}\rho \\
&\le e^{Cn\la_H(A)+A}\rho(\bT) \qquad \qquad
\bigl( {\rm by\ extremality\ of\ }
P, \int_{\bT} |P|^2\, {\rm d}\rho \le \rho(\bT)\, \bigr)
\\
&\le e^{CA} \quad\qquad\qquad\qquad\qquad ({\rm by\ \eqref{ko1}\ and\ the\ equality\ }n\la_H(A)=A),
\end{align*}
whence,
\[
\la_H(A)\, \log\Bigl( \int_\bT |P|^2\, {\rm d}m \Bigr)
\le C A\la_H(A).
\]
Therefore,
\[
0 \le \frac2e + C A\la_H(A) + 2\log e_n(\rho) + \int_{\{H\le A\}} H\, {\rm d}m,
\]
and finally,
\[
| \log e_n(\rho)| \lesssim A\la_H(A) + \int_{\{H\le A\}} H\, {\rm d}m
= \int_\bT H_A\, {\rm d}m,
\]
proving Theorem~\ref{thm:LB}. \hfill $\Box$

\section{The upper bound for $e_n$ via Taylor polynomials of an outer function}

We give two upper bounds for $e_n(\rho)$. Both of them are based on the construction
of monic polynomials of large degree with a good estimate for the $L^2(\rho)$-norm.
The first bound uses Taylor polynomials of an outer
function $F$ such that $1/F$ mimics the behaviour of $\rho$. It is better adjusted to
the case when the distribution function $\la_H(a)$ decays relatively fast as $a\to\infty$.
The second bound uses classical Chebyshev's polynomials and starts working
only when $\la_H (a)$ decays at infinity slower than $1/a$.

\medskip
Let $\phi\colon [0, \frac12]\to (0, +\infty]$ be a continuous decreasing function, $\phi (0)=\infty$, $\phi(\tfrac12)\le \inf_{\bT} H$.
Given $\tau\in [-\frac12, \frac12]$, denote by $\theta_\tau$ solution to the equation
$\phi(\theta_\tau) = H(e^{2\pi{\rm i}\tau})$.
We call the function $H$ {\em subordinated to} $\phi$ if, for any $\tau$,
\begin{align*}
H(e^{2\pi{\rm i}\theta}) &\le \phi(\theta+\theta_\tau-\tau),  \qquad \tau-\theta_\tau < \theta \le \tau, \\
H(e^{2\pi{\rm i}\theta}) &\le \phi(\tau+\theta_\tau-\theta),  \qquad \tau \le \theta < \tau+\theta_\tau.
\end{align*}
Note that an equivalent way to express the $\phi$-subordination is to say that the
function
$(\phi^{-1}\circ H)(e^{2\pi{\rm i}\theta})$ is a non-negative
Lipschitz function on $[-\frac12, \frac12]$ with the Lipschitz constant at most one.

We call the unbounded continuous decreasing function $\phi$ on $(0, \frac12]$ {\em regular} if it satisfies at least one of the following two conditions:
\[
{\rm the\ function\ }
\theta\mapsto \theta\phi(\theta)
\ {\rm does\ not\ decrease
\ and\ } \phi(\theta)\gtrsim \log\frac1{\theta}
\eqno({\rm Reg}1)
\]
\[
\phi(\theta/2)\lesssim \phi(\theta) \ {\rm and\ }
\phi(\theta) \gtrsim 1/\theta.
\eqno({\rm Reg}2)
\]

\begin{theorem}\label{thm:UB1}
Suppose that
\[
\int_\bT e^H\, {\rm d}\rho <\infty,
\]
with $H$ subordinated to a regular function $\phi$. Then
\[
| \log e_n (\rho) | \gtrsim \int_0^A \la_H(a)\,{\rm d}a,
\]
where $A$ solves the equation $n\phi ^{-1}(A)=A$ when $\phi$ satisfies condition {\rm (Reg1)},
and $A=\sqrt{n}$ when $\phi$ satisfies condition {\rm (Reg2)}.
\end{theorem}

\begin{corollary}\label{cor:UB1} In the assumptions of Theorem~\ref{thm:UB1}, suppose that 
$\log \phi(\theta)\gtrsim \log\frac1{\theta}$.
If
\[
\liminf_{a\to\infty} a\la_H(a) > 0,
\]
then $e_n(\rho)$ decay to zero at least as a negative power of $n$.

Furthermore, $e_n(\rho)$ decay to zero faster than any negative power of $n$, provided
that
\[
\lim_{a\to\infty} a\la_H(a) = \infty.
\]
\end{corollary}

Note that we need to impose the additional condition $\log\phi(\theta) \gtrsim \log\tfrac1{\theta}$ in this Corollary 
only in the case when $\phi$ satisfies the first regularity condition (Reg1).

\subsection{Taylor polynomials}
Denote by $\bP_r$ the Poisson kernel for the unit disk evaluated at the point $r\in [0, 1)$.

\begin{lemma}\label{lemma_Taylor}
Let $H$ be a weight such that
\[
H_A * \bP_{1-\delta} \lesssim H + 1 \qquad {\rm everywhere\ on\ } \bT,
\]
with
$ \log \delta^{-1} \lesssim A $.
Suppose that
\[
\int_\bT e^H\, {\rm d}\rho < \infty.
\]
Then there exists a positive constant $C$
such that, for $n\ge CA/\delta$,
we have
\[
|\log e_n(\rho) | \gtrsim \int_\bT H_A\, {\rm d}m.
\]
\end{lemma}

\subsubsection{Proof of Lemma~\ref{lemma_Taylor}}

Let $M$ be a positive constant such that
$ H_A * \bP_{1-\delta} \le M(H + 1)$, and
let $F_A$ be an outer function in $\bD$ with the boundary values
$|F_A|^2 = e^{H_A/M}$, i.e.,
\[
\log |F_A(rt)| = \frac1{2M} \bigl( H_A * \bP_r \bigr) (t).
\]
We expand $F_A((1-\delta)z)$ into the Taylor series
\[
F_A((1-\delta)z) = \sum_{k\ge 0} f_k z^k,
\]
and consider the Taylor polynomials
\[
P_A(z) = \sum_{k=0}^n f_k z^k.
\]
Then,
\[
e_n(\rho)^2 \le |P_A(0)|^{-2}\, \int_\bT |P_A|^2\, {\rm d}\rho.
\]

First, we note that
\[
|P_A(0)| = |F_A(0)| = \exp\Bigl( \int_\bT \log |F_A|\, {\rm d}m \Bigr)
= \exp\Bigl( \frac1{2M}\, \int_\bT H_A \, {\rm d}m \Bigr),
\]
i.e.,
\[
e_n(\rho)^2 \le \int_\bT |P_A|^2\, {\rm d}\rho \cdot
\exp\Bigl( - \frac1{M}\, \int_\bT H_A \, {\rm d}m \Bigr).
\]
Next,
\[
\int_\bT |F_A((1-\delta)t)|^2\, {\rm d}\rho(t)
= \int_\bT \exp \Bigl( \frac1{M} H_A * \bP_{1-\delta} \Bigr)\, {\rm d}\rho
\le \int_\bT \exp\bigl( H + 1\bigr)\, {\rm d}\rho
\lesssim 1,
\]
so it remains to estimate the remainder
\[
| F_A((1-\delta)z) - P_A(z) | \le \sum_{k>n} |f_k|.
\]
By Cauchy's estimates,
\[
|f_k| \le (1-\delta)^{k} \max_{\bT} \, |F_A|
\le (1-\delta)^{k} e^{A/(2M)},
\]
whence,
\[
\sum_{k>n} |f_k| \lesssim \delta^{-1} e^{A/(2M)-n\delta} \lesssim 1\,,
\]
provided that $A/\delta \lesssim n$ (here, we use that $\log \delta^{-1} \lesssim A$).
Thus,
\[
\int_\bT |P_A|^2\, {\rm d}\rho \lesssim 1,
\]
which proves the lemma. \hfill $\Box$

\subsection{Estimates of the Poisson integral}

Put $p_\delta (\theta) = \bP_{1-\delta}(e^{2\pi{\rm i}\theta})$ and recall that
$p_\delta(\theta) \lesssim \min(\delta^{-1}, \delta\theta^{-2})$.

\begin{lemma}\label{lemma_Poisson_A}
Let $\phi\colon (0, \frac12]\to [0, \infty)$ be an unbounded continuous decreasing function, let
$\widetilde\phi$ be its even $1$-periodic extension on $\bR$, and
$\widetilde\phi_A=\min(\widetilde\phi, A)$. Then
\[
\widetilde\phi_A * p_\delta \lesssim \widetilde\phi + 1
\qquad {\rm everywhere\ on\ } \bR,
\]
provided that at least one of the following holds:

\noindent$\rm (i)$ the function $\theta\mapsto \theta^2\phi(\theta)$ does not decrease, and
$\delta\lesssim\phi^{-1}(A)$;

\noindent$\rm (ii)$
$\phi(\theta/2)\lesssim \phi(\theta)$, $\phi(\theta)\gtrsim 1/\theta$,
and $\delta\lesssim 1/A$.
\end{lemma}
Note that condition (i) is weaker than condition (Reg1) in Theorem~\ref{thm:UB1},
i.e., the lemma is a bit stronger than what we will use for the proof of Theorem~\ref{thm:UB1}. We need this version of Lemma~\ref{lemma_Poisson_A} for
the proof of Theorem~\ref{thm:DeepZero}.
We also note that condition (i) yields estimate $\phi(\theta/2)\lesssim\phi(\theta)$
from condition (ii).

\subsubsection{Proof of Lemma~\ref{lemma_Poisson_A}}
We take a sufficiently small $\tau_0>0$ so that $\phi(\tau_0)\ge 1$,
fix $\tau\in (0, \tau_0]$, and estimate the convolution
$(\widetilde\phi_A *p_\delta)(\tau)$.
There is nothing to prove if $A\le\phi\bigl( \frac12 \tau\bigr)$ since in this case
\[
(\widetilde\phi_A * p_\delta)(\tau) \le \max_{[-\frac12, \frac12]} \widetilde\phi_A
= A \le
\phi\bigl( \tfrac12 \tau \bigr) \lesssim
\phi(\tau)
\]
for any $\delta>0$. Hence, in what follows, we assume that $A\ge\phi\bigl( \frac12 \tau\bigr)$, i.e., $\tau\ge 2\phi^{-1}(A)$.

First, we note that for any $\theta\in [0, \frac12]$, we have
$\widetilde\phi_A(\tau+\theta) \le \widetilde\phi_A(\tau-\theta)$
(to see this, one needs to consider three cases: $0\le\theta\le\tau$,
$\tau\le\theta\le\frac12-\tau$, and $\frac12-\tau\le\theta\le\frac12$).
Therefore,
\begin{multline*}
\int_{-\frac12}^{\frac12} \widetilde\phi_A(\tau-\theta) p_\delta(\theta)\, {\rm d}\theta
= \int_0^{\frac12}
\bigl[
\widetilde\phi_A(\tau-\theta) + \widetilde\phi_A(\tau+\theta)
\bigr]
p_\delta(\theta) \, {\rm d}\theta \\
\le 2A\, \int_{|\tau-\theta|\le \phi^{-1}(A)} p_\delta(\theta)\, {\rm d}\theta
+ 2\, \int_
{\substack{|\theta-\tau|\ge \phi^{-1}(A), \\
0\le\theta\le\frac12}}
\phi(\tau-\theta) p_\delta(\theta)\, {\rm d}\theta
= I+ II.
\end{multline*}
Before we start estimating integrals on the RHS,
observe that $\delta\lesssim\tau$. In the case (i) it is obvious since
$\delta\lesssim \phi^{-1}(A)\le\tau/2$, in the case (ii) it is also obvious since
then $\delta\lesssim 1/A \lesssim \phi^{-1}(A) \le \tau/2$.
Therefore, in the first integral
$\theta \ge \tau-\phi^{-1}(A) \ge \tau/2 \gtrsim \delta$.
Recalling the standard estimate of the Poisson kernel
$p_\delta(\theta) \lesssim \min (\delta^{-1}, \delta \theta^{-2}) $,
we get
\[
I \lesssim
A\delta\, \int_{\tau-\phi^{-1}(A)}^{\tau+\phi^{-1}(A)} \frac{{\rm d}\theta}{\theta^2}
\lesssim \frac{A\delta\phi^{-1}(A)}{\tau^2}.
\]
In both cases (i) and (ii) the RHS is $\lesssim\phi(\tau)$. Indeed, if (i) holds, then it
is bounded by $ \phi^{-1}(A)^2 A/\tau^2 \le \tau^2\phi(\tau)/\tau^2 = \phi(\tau) $.
If (ii) holds, then it is bounded by
$ \phi^{-1}(A)/\tau^2 \lesssim 1/\tau \lesssim \phi(\tau) $.

We split the second integral into four parts
\[
\int_{\substack{|\theta-\tau|\ge \phi^{-1}(A), \\
0\le\theta\le\frac12}} =
\int_0^{\min(\delta, \frac12 \tau)} +
\int_{\min(\delta, \frac12 \tau)}^{\frac12 \tau} +
\int_{\phi^{-1}(A)\le |\theta-\tau| \le \frac12 \tau} +
\int_{\frac32 \tau}^{1/2}
\]
and estimate them one by one. We have
\begin{multline*}
\int_0^{\min(\delta, \frac12 \tau)} \phi(\tau-\theta) p_\delta(\theta)\, {\rm d}\theta
\lesssim \frac1{\delta}\, \int_0^{\min(\delta, \frac12 \tau)}
\phi(\tau-\theta)\, {\rm d}\theta \\
\le \frac1{\delta}\, \min(\delta, \tfrac12 \tau) \cdot \phi(\tfrac12 \tau)
\lesssim \phi(\tau),
\end{multline*}
and
\begin{multline*}
\int_{\min(\delta, \frac12 \tau)}^{\frac12 \tau} \phi(\tau-\theta) p_\delta(\theta)\,
{\rm d}\theta \lesssim
\delta\, \int_{\min(\delta, \frac12 \tau)}^{\frac12 \tau} \frac{\phi(\tau-\theta)}{\theta^2}
\, {\rm d}\theta \\
\lesssim \delta \cdot \phi(\tfrac12 \tau)\, \int_{\min(\delta, \frac12 \tau)}^\infty
\frac{{\rm d}\theta}{\theta^2}\,  \stackrel{\delta\lesssim\tau}\lesssim \, \phi(\tau).
\end{multline*}
Next,
\begin{multline*}
\int_{\phi^{-1}(A)\le |\theta-\tau| \le \frac12 \tau}
\phi(\tau-\theta) p_\delta(\theta)\, {\rm d}\theta \\
\lesssim
\delta\, \int_{\phi^{-1}(A)\le |\theta-\tau| \le \frac12 \tau}
\frac{\phi(\tau-\theta)}{\theta^2}\, {\rm d}\theta
\simeq \frac{\delta}{\tau^2}\, \int_{\phi^{-1}(A)}^{\frac12 \tau}
\phi(\xi)\, {\rm d}\xi.
\end{multline*}
In the case (i), the integral on the RHS equals
\[
\frac{\delta}{\tau^2}\, \int_{\phi^{-1}(A)}^{\frac12 \tau}
\frac{\xi^2 \phi(\xi)}{\xi^2}\, {\rm d}\xi
\le \frac{\delta}{\tau^2}\, \frac{(\tau/2)^2 \phi(\tau/2)}{\phi^{-1}(A)}
\, \stackrel{\delta\lesssim \phi^{-1}(A)}\lesssim\, \phi(\tau),
\]
while in the case (ii), it does not exceed
\[
\frac{\delta}{\tau^2} \cdot A\tau/2 \,
\stackrel{\delta\lesssim 1/A}\lesssim \, \frac1{\tau}
\, \stackrel{\tau\phi(\tau)\gtrsim 1}\lesssim \, \phi(\tau).
\]
At last,
\[
\int_{\frac32 \tau}^{1/2} \phi(\tau-\theta) p_\delta(\theta)\, {\rm d}\theta
\lesssim \delta \int_{\frac32 \tau}^{\frac12} \frac{\phi(\theta - \tau)}{\theta^2}\,
{\rm d}\theta
\lesssim \delta \phi(\tau/2) \cdot \frac1\tau \,
\stackrel{\delta\lesssim\tau}\lesssim\, \phi(\tau),
\]
completing the proof of the lemma. \hfill $\Box$

\subsubsection{The Poisson integral of $H_A$}

\begin{lemma}\label{lemma:Poisson} Let $H$ be subordinated to a regular function $\phi$,
and $H_A=\min (H, A)$. Then
\[
H_A * \bP_{1-\delta} \lesssim H+1 \quad{\rm everywhere\ on\ } \bT,
\]
provided that $\delta\lesssim\phi^{-1}(A)$ when $\phi$ satisfies condition {\rm (Reg1)},
and $\delta\lesssim A^{-1}$ when $\phi$ satisfies condition {\rm (Reg2)}.
\end{lemma}

Clearly, Lemma~\ref{lemma_Taylor} and Lemma~\ref{lemma:Poisson} combined together
yield Theorem~\ref{thm:UB1}.

\paragraph{Proof of Lemma~\ref{lemma:Poisson}}
We write $H(e^{2\pi{\rm i}\theta})=h(\theta)$, fix the point
$\tau\in [-\frac12, \frac12]$ with $h(\tau)<\infty$ at which we will
estimate the convolution $(h_A * p_\delta)(\tau)$,
and choose $\theta_\tau$ so that $\phi(\theta_\tau)=h(\tau)$. Similarly to the proof of
the previous lemma, we assume that $A \ge \phi(\frac12 \theta_\tau)$, i.e.,
that $\phi^{-1}(A) \le \frac12 \theta_\tau$; otherwise,
\[
(h_A * p_\delta)(\tau) \le A \le \phi(\tfrac12 \theta_\tau)
\lesssim  \phi(\theta_\tau) = h(\tau),
\]
and we are done. Now,
\[
(h_A * p_\delta)(\tau) =
\Bigl(\, \int_{|\theta-\tau|\ge \theta_\tau-\phi^{-1}(A)} + \,
\int_{|\theta-\tau|\le \theta_\tau-\phi^{-1}(A)} \,
\Bigr)\, h_A(\theta) p_\delta(\tau-\theta)\, {\rm d}\theta = I + II.
\]

To estimate the first integral, we note that, since
$\theta_\tau-\phi^{-1}(A)\ge \theta_\tau/2$, we have
\[
I \le A \int_{|\theta-\tau|\ge \theta_\tau/2} \, p_\delta(\tau-\theta)\, {\rm d}\theta
\le 2A \int_{\theta_\tau/2}^{1/2} p_\delta (\theta)\, {\rm d}\theta
\lesssim A\delta \int_{\theta_\tau/2}^\infty \frac{{\rm d}\theta}{\theta^2}
\lesssim \frac{A\delta}{\theta_\tau}\,.
\]
In the first case, the RHS is
\[
\, \stackrel{\delta\lesssim \phi^{-1}(A)}\lesssim\,
\frac{A\phi^{-1}(A)}{\theta_\tau}
\, \stackrel{\phi^{-1}(A)\le \frac12 \theta_\tau}\lesssim \,
\frac{\frac12 \theta_\tau \phi(\frac12 \theta_\tau)}{\theta_\tau}
\lesssim \phi(\theta_\tau).
\]
In the second case, $A\delta/\theta_\tau \lesssim 1/\theta_\tau
\lesssim \phi(\theta_\tau)$. Therefore, in both cases, the first integral is
$\lesssim \phi(\theta_\tau)=h(\tau)$.

To estimate the second integral, we note that, by the subordination to $\phi$,
it is bounded by
$ 2 (\widetilde\phi_A * p_\delta)(\theta_\tau) $, which, by the previous lemma, is
$\lesssim \phi(\theta_\tau)+1 = h(\tau)+1$.
\hfill $\Box$

\section{The upper bound for $e_n$ via Chebyshev polynomials}
\label{sec5}

Here we assume that the function $H$ is lower semicontinuous; i.e., the sets
$\{H\!>\!a\}$ are open, and denote by $\la^*_H(a)$ the length of the longest open interval
within $\{H\!>\!a\}$.

\begin{theorem}\label{thm:UB2}\label{the9}
Suppose that
\[
\int_\bT e^H\, {\rm d}\rho <\infty.
\]
Then
\[
| \log e_n (\rho) | \gtrsim A \la_H^*(A)\,,
\]
where $A=A(n)$ solves the equation $n\la_H^*(A)=A$.
\end{theorem}

The following Corollary combines Theorem~\ref{thm:UB2} with Theorem~\ref{thm:LB}
(and takes into account Observations~2.1 and~2.3)
\begin{corollary}\label{cor:UB2}
Let ${\rm d}\rho = e^{-H}\, {\rm d}m $. Suppose
that the set $\{H>a\}$ contains an interval with the length comparable to the
total length of $\{H>a\}$ (i.e. $\la_H^* \gtrsim \la_H$), and that the function $\la_H$ satisfies
\[
\limsup_{a\to\infty} \frac{\la_H(a)}{\la_H(2a)} <2.
\]
Then
\[
| \log e_n (\rho)| \simeq A\lambda_H(A)\,,
\]
where $A=A(n)$ is a solution to the equation $n\la_H(a)=a$.
\end{corollary}

\subsection{Proof of Theorem~\ref{thm:UB2}}

We will use the following classical lemma (cf., for instance,~\cite{MR}):
\begin{lemma}\label{lemma:Chebyshev-Maergoiz}
For any $\alpha\in (0, \frac{\pi}2 ]$ and for any $n\in 2\bN$,
there exists a monic polynomial $T_{n,\alpha}$ of degree $n$ such that
$$
\max\{|T_{n,\alpha}(e^{i\theta})|\colon |\theta| \ge \alpha \} = 2\cos^n(\alpha/2)\,.
$$
\end{lemma}

\noindent For the reader's convenience, we recall its proof. Put
$$
T_{2m,\alpha}(e^{i\theta}) = 2\cos^{2m}(\alpha/2)e^{im\theta}
\cos\left( 2m\arccos\left(\frac{\cos(\theta/2)}{\cos(\alpha/2)} \right)\right)\,.
$$
We only need to show that this is a monic polynomial of degree $n=2m$.
Recall that $\cos(2m\arccos x) = 2^{2m-1}Q_m(x^2)$, where $Q_m$ is a monic polynomial
of degree $m$. Then
\begin{align*}
T_{2m,\alpha}(e^{i\theta}) &= 2^{2m} \cos^{2m}(\alpha/2)e^{im\theta}
Q_m\left(\frac{\cos^2\theta/2}{\cos^2\alpha/2}\right) \\
&=
2^{2m}\cos^{2m}(\alpha/2)e^{im\theta}
Q_m\left(\frac{e^{i\theta}+e^{-i\theta}+2}{4\cos^2\alpha/2}\right)
\end{align*}
and it is easily seen  that the RHS is a monic polynomial of degree $2m$,
proving the lemma. \hfill $\Box$

\medskip
Now, we turn to the proof of Theorem~\ref{thm:UB2}. Without loss of generality, we assume that $n$ is
an even number.
Let $J\subset\bT$ be the longest arc in the set $\{H>a\}$. We assume that
$J=\{t=e^{{\rm i}\varphi}\colon |\varphi|\le\alpha \}$, $\alpha = \pi\la_H^*(a)$.
Let $T=T_{n, \alpha}$ be a monic polynomial
of degree $n$ as in Lemma~\ref{lemma:Chebyshev-Maergoiz}. Then, by a straightforward computation
(or by the classical Remez inequality)
$$
\max_J |T| = \max_{|\varphi|\le \alpha} \left| T(e^{i\varphi})\right|
\lesssim  e^{Cn\alpha} \, \bigl( \cos\frac{\alpha}{2} \bigr)^n\,.
$$
Noting that $ \cos^n\tfrac{\alpha}{2} \le e^{-cn\alpha^2}$, we get
\begin{multline*}
\int_\bT |T|^2\, {\rm d}\rho \le \max_{\bT} \bigl( |T|^2e^{-H} \bigr)\,
\int_\bT e^H\, {\rm d}\rho \\
\lesssim e^{-a}\, \max_J |T|^2 + \max_{\bT\setminus J} |T|^2
\lesssim \Bigl( e^{Cn\la_H^*(a) - a} + 1 \Bigr) e^{-cn\la_H^*(a)^2}.
\end{multline*}
Letting $A_C$ be the unique solution to the equation $Cn \la_H^*(A_C)=A_C$,
we obtain
\[
| \log e_n(\rho) | \gtrsim
n \la_H^*(A_C)^2
\simeq A_C\la_H^*(A_C)
\simeq A\la_H^*(A),
\]
completing the proof of Theorem~\ref{thm:UB2}. \hfill $\Box$

\section{Examples}

To illustrate our results, we consider the function $H=h\circ d_K$, where
$h\colon (0, \frac12]\to (0, +\infty)$ is a $C^1$-smooth decreasing function, $h(0)=+\infty$,
and $d_K(t)=\operatorname{dist}(t, K)$, where $K\subset \bT$ is a compact set of zero length
(recall that we identify $\bT$ with $\bR/\bZ$).

Denote by $K_{+s}=\{t\colon d_K(t)<s \}$ the $s$-neighbourhood of $K$ and by
$\psi_K(s) = m(K_{+s})$ its length. Then
\[
\la_H(a) = m\{H>a\} =
\begin{cases}
\psi_K(h^{-1}(a)), &a\ge h(\tfrac12) \\
1, &a<h(\tfrac12),
\end{cases}
\]
and
\[
\int_0^A \la_H(a)\, {\rm d}a= \Bigl(\int_0^{h(1/2)}+\int_{h(1/2)}^A\Bigr) \la_H(a)\, {\rm d}a= \int_{h^{-1}(A)}^{1/2} \psi_K |h'| + h(\tfrac12),
\]
provided that $A>h(\frac12)$.

To estimate the function $\psi_K$ it is convenient to use that
$ \psi_K(s) \simeq s N_K(s) \simeq sP_K(s) $,
where $N_K(s)$ is the covering  number of $K$ and $P_K(s)$ is the packing number
of $K$, see, for instance,~\cite[Chapter~3]{Falc}.
We call the set $K$ $\gamma$-regular if
$\psi_K(s) \simeq s^{1-\gamma}$. For instance,
the set $e^{2\pi {\rm i}\mathcal C}$, where
$\mathcal C$ is the standard ternary Cantor set is
$\gamma$-regular with
$\gamma = \frac{\log 2}{\log 3}$,
while the set
$\{t=\exp(2\pi {\rm i} n^{-\nu})\colon n\in\bN \} \cup \{1\}$
is $\gamma$-regular with $\gamma=(\nu+1)^{-1}$.

\subsection{Two corollaries}
We get straightforward corollaries to our results taking $h(s)=s^{-p}$.
\begin{corollary}\label{cor:d_K2}
Let $K\subset\bT$ be a $\gamma$-regular compact set with some $\gamma\in [0, 1)$.
Suppose that
${\rm d}\rho = \exp \bigl( - d_K^{\gamma-1} \bigr)\, {\rm d}m$. Then
$ | \log e_n(\rho)| \simeq \log n$.
\end{corollary}

The second corollary pertains to the case when
the length of the longest interval in the set $\{d_K<s\}$ is comparable with the
length of the whole set $\{d_K<s\}$. Then Corollary~\ref{cor:UB2}
applies.
\begin{corollary}\label{cor:d_K3}
Let $\nu>0$, $K=\{t=\exp(2\pi {\rm i}n^{-\nu})\colon n\in\bN \} \cup \{1\}$, and
${\rm d}\rho = \exp \bigl( - d_K^{-p} \bigr)\, {\rm d}m$ with
$p>\vartheta = \frac{\nu}{\nu+1}$. Then
\[
|\log e_n (\rho)| \simeq n^{\frac{p-\vartheta}{p+\vartheta}}.
\]
\end{corollary}

\subsection{Measures with deep zero at one point}
The last illustration to our estimates  pertains to the simplest
case when the measure $\rho$ has a deep zero at one point and is symmetric
with respect to this point. In this case, our estimates yield a relatively
complete result.
\begin{theorem}\label{thm:DeepZero} Let $h\colon (0, \frac12]\to [0, +\infty)$ be a
continuous decreasing function such that
\[
\int_0 h(a)\,{\rm d}a = +\infty.
\]
Suppose that $h$ satisfies at least one of the following two conditions:
\begin{itemize}
\item[{\rm (i)}]
\[
\theta\mapsto \theta^2 h(\theta) \quad {\rm does\ not\ decrease},
\]
and
\[
|\log \theta| = O(h(\theta)),\qquad \theta\to 0;
\]
\item[{\rm (ii)}]
\[
\limsup_{a\to\infty} \frac{h^{-1}(a)}{h^{-1}(2a)} < 2
\]
\end{itemize}
Let $\rho$ be an absolutely continuous measure on $\bT$ with density
$e^{-h(|\theta|)}$.
Then
\[
|\log e_n (\rho)| \simeq \int_0^{1/2} h_A(a)\,{\rm d}a
\]
where $A$ solves the equation 
$n h^{-1}(A)=A$ and $h_A=\min(h, A)$.
\end{theorem}

The lower bound for $e_n$ (i.e., the upper bound for $|\log e_n|$) follows from Theorem~\ref{thm:LB} and does not need any regularity assumptions on $h$.
Conditions (i) and (ii) are needed for the proof of the upper bound for $e_n$.
In the case (i), it is a consequence of Lemma~\ref{lemma_Taylor} combined with
the first case of Lemma~\ref{lemma_Poisson_A}. In the case~(ii), it follows from Theorem~\ref{thm:UB2}. Note that these two cases overlap, e.g.,
the function $h(\theta)=\theta^{-p}$ with $1<p\le 2$ satisfies both of them.

The following corollary gives an idea about the rate of decay of $e_n(\rho)$
for several explicitly written functions $h$.

\begin{corollary} Let $\rho$ be an absolutely continuous measure on $\bT$ with density
$e^{-h(|\theta|)}$. Then for $n\ge 4$ we have
\begin{itemize}
\item[{\rm (i)}] If $h(\theta)\simeq \theta^{-1}\log^{-1}(1/\theta)$, then $|\log e_n (\rho)| \simeq \log\log n$,
\item[{\rm (ii)}] If $p>-1$ and $h(\theta)\simeq \theta^{-1}\log^p(1/\theta)$, then $|\log e_n (\rho)| \simeq (\log n)^{p+1}$,
\item[{\rm (iii)}] If $p>1$ and $h(\theta)\simeq \theta^{-p}$, then $|\log e_n (\rho)| \simeq n^{(p-1)/(p+1)}$,
\item[{\rm (iv)}] If $p>0$ and $h(\theta)\simeq \exp(\theta^{-p})$, then $|\log e_n (\rho)| \simeq n(\log n)^{-2/p}$.
\end{itemize}
\end{corollary}

\medskip
Our last remark is that, plausibly, the technique based on the potential theory
in the external field developed by Mhaskar--Saff, Rakhmanov, Levin--Lubinsky, Totik
and others should allow one to obtain more precise estimates of $e_n$ in the situation
considered in Theorem~\ref{thm:DeepZero}. See for instance, Theorem~1.22 and Examples~3 and~4 in Section~1.6 in \cite{LL} which contain similar results for orthogonal polynomials on the real line. On the other
hand, likely, this will require much stronger regularity assumptions
on the function $h$ and more technical proofs.

\bigskip
\medskip

\noindent{A.B.: Institut de Math\'ematiques de Marseille,
Aix Marseille Universit\'e, CNRS, Centrale Marseille, I2M, Marseille, France
\newline {\tt alexander.borichev@math.cnrs.fr}
\smallskip
\newline\noindent A.K.: Department of Mathematics and Mechanics,
St. Petersburg State University, St. Petersburg, Russia
\newline {\tt a.kononova@spbu.ru}
\smallskip
\newline\noindent M.S.:
School of Mathematics, Tel Aviv University, Tel Aviv, Israel
\newline {\tt sodin@tauex.tau.ac.il}
}
\end{document}